\documentclass[12pt]{amsart}
\usepackage{amssymb,latexsym, amsmath, amscd, array, graphicx}
\usepackage{amssymb}

\usepackage{setspace}
\input xy
\xyoption{all}
\makeatletter 

\usepackage{setspace}
\usepackage{enumitem}

\date{}

\newtheorem{theorem}{Theorem}[section]
\newtheorem{claim}[theorem]{Claim}

\newtheorem{proposition}[theorem]{Proposition}
\newtheorem{definition}[theorem]{Definition}

\newcommand{\U}{{\Bbb U}}

\newcommand{\B}{{\Bbb B}}

\newcommand{\N}{{\Bbb N}}

\newcommand{\w}{{\rm w}}

\newcommand{\st}{{\rm st}}

\newcommand{\lo}{\rightarrow}

\newcommand{\black}{{\blacksquare}}

\newcommand{\diam}{{\rm diam}}

\newcommand{\asdim}{{\rm asdim}}

\newcommand\blfootnote[1]{%
  \begingroup
  \renewcommand\thefootnote{}\footnote{#1}%
  \addtocounter{footnote}{-1}%
  \endgroup
}
\begin{document}

\title{\bf Mardesic factorization theorem for asymptotic dimension }

\author{ Jerzy Dydak, Michael Levin and Jeremy Siegert}
\blfootnote{The second and third authors
 were supported by the Israel Science Foundation grant No. 2196/20 and 
 the Institute of Mathematics of the Polish Academy of Sciences where 
 the main results were obtained during their visits in 2021-2022.}

\maketitle

\begin{abstract} 
The main goal of this note is to prove a coarse analogue  of Factorization Theorems in Dimension Theory:\\
\emph{Let $f: X \lo Y$ be a coarsely continuous map. Then $f$ factors through coarsely continuous maps
$g : X \lo Z$ and $h : Z \lo Y$ with asymptotic dimension of $Z$ at most asymptotic dimension of $X$ and  the weight of  $Z$ at most the weight of $Y$.}
\end{abstract}
\begin{section}{Introduction}  

We say that a map $f : X \lo Y$  between metric spaces 
is \textbf{coarsely continuous} if $f$ sends uniformly  bounded covers in $X$
to uniformly bounded  families in $Y$.  Maps $f,g : X \lo Y$ are \textbf{coarsely close} if 
there is $r>0$ such that
   $f(x)$ and $g(x)$  are $r$-close for all  $x \in X$.
A coarsely continuous  map $ f : X\lo Y$ is a \textbf{coarse equivalence} if there is 
a coarsely continuous map $g : Y \lo X$  such that $g\circ f$ and $f\circ g$ are coarsely close
to the identity maps of $X$ and $Y$ respectively. A map  $f : X \lo Y$ is a \textbf{coarse embedding} if $f$ 
is a coarse equivalence between $X$ and $f(X)$. Note that a coarsely continuous map
$f: X\lo Y$ is a coarse embedding if and only if the preimages of uniformly bounded covers of $Y$ are
uniformly bounded in $X$.
 We recall that  the \textbf{asymptotic dimension}    of a metric space $X$  is bounded by $n$, written $\asdim X \leq n$,
  if for every
$r>0$ there is a uniformly bounded cover of $X$ that splits into $n+1$ families of $r$-disjoint sets.
The weight of $X$
is denoted by $\w X$. All the spaces are assumed to be metric.

The main goal of this note is to prove 
\begin{theorem}
\label{factorization}
Let $f: X \lo Y$ be a coarsely continuous map. Then $f$ factors through coarsely continuous maps
$g : X \lo Z$ and $h : Z \lo Y$ with $\asdim Z \leq \asdim X$ and   $\w Z \leq \w Y$.
\end{theorem}

Theorem \ref{factorization} is a coarse analogue  of factorization theorems in Dimension Theory 
which originated in Mardesic's work \cite{mardecis}.  Using a well-known approach we derive from
Theorem \ref{factorization}

\begin{theorem}
\label{universal}
{\rm (Bell-Nag\' orko \cite{bell-nagorko})} For every $n$ there is a separable metric space 
of $\asdim=n$ which is universal for separable metric spaces of $\asdim \leq n$ (i.e. it contains
a coarsely equivalent copy of every separable metric space of $\asdim \leq n$).
\end{theorem}

Note that universal  spaces  for $0$-dimensional  spaces of bounded geometry and 
$0$-dimensional proper spaces were constructed 
in \cite{duniversal} and \cite{dydakuniversal}.

This note is organized as follows:  in the next section we  present auxiliary results, Theorem \ref{universal}
 is proved in section 3 and   Theorem \ref{factorization} is proved in Section 4.

\end{section}
\begin{section}{Witnessing and defining asymptotic dimension}

Let  $\mathcal{F}$  be a family of subsets  of $X$.  By  the \textbf{star} $\st  (A, \mathcal{F})$   of a subset $A$ of  $X$  
we denote
the union of $A$ with all the sets $F \in \mathcal{F}$ meeting  $A$,  and $\st \mathcal{F}$ stands for the cover of $X$
consisting of $\st (x, {\mathcal{F}})$ for all $ x \in X$. We say that $\mathcal{F}$ \textbf{separates} subsets 
$A$ and $B$ of $X$ if $\st(A, {\mathcal{F}})$ does not meet $B$.

For a metric space $(X,d)$ and  a subset $A\subset X$ we denote by $\B(A,r)$ the set of  points 
at distance at most $r$ from $A$.   We say that a subset $B$ of $X$ is \textbf{$r$-close} to $A$ if 
in $B \subset \B(A,r)$. 

\begin{proposition}
\label{prop-for-witnessing}
Let $X$ be a metric space with $\asdim X \leq n$. Then for every $r>0$ there is 
a uniformly bounded cover  of $X$ which has a  Lebesgue number 
 $r$ and  which   splits into $n+1$ families  of $r$-disjoint sets.
Moreover if $b$ is  a base point in $X$ then we can assume that $\B(b,r)$ is contained
in an element of each family of this splitting.
\end{proposition}
{\bf Proof.} Let $\mathcal{F}$ be a uniformly bounded cover of $X$ that splits
into families ${\mathcal{F}}_j, 1\leq j \leq n+1,$ of $(3r)$-disjoint sets.
Then the families ${\mathcal{F}}'_j=\{\B(F,r): F \in {\mathcal{F}}_j \}$ are $r$-disjoint and form
 the cover ${\mathcal{F}}'={\mathcal{F}}'_1 \cup \dots\cup {\mathcal{F}}'_{n+1}$ whose
 Lebesgue number is  $r$.

 Now suppose  $b$ is a base point of $X$. Define
  ${\mathcal{F}}''_j$ to be the family consisting of elements of ${\mathcal{F}}'_j $ not intersecting $B(b,2r)$ 
  and one extra element being the the union of $B(b,r)$ with the elements of ${\mathcal{F}}'_j$ intersecting $B(b, 2r)$.
   Then 
 the cover ${\mathcal{F}}''={\mathcal{F}}''_1 \cup \dots\cup {\mathcal{F}}''_{n+1}$ satisfies
  the second conclusion of the proposition. $\black$

\begin{definition}
Let $X$ be  a metric space of  $\asdim\leq n$.  A sequence $\mathcal{F}$ of covers  ${\mathcal{F}}_i, i\in \N,$ 
of $X$ is said to 
\textbf{witness} $\asdim \leq n$ if   ${\mathcal{F}}_{i}$ is $R_{i}$-bounded,  $r_i$ is
a Lebesgue number of ${\mathcal{F}}_{i}$
 and ${\mathcal{F}}_{i}$ splits into the union of 
$n+1$ families ${\mathcal{F}}_{ij}, 1\leq j \leq n+1$, of $ r_i$-disjoint  sets  
with $ r_{i+1}> (100i +1)R_{i}$
and $R_i  > r_i$.
We will say that $\mathcal{F}$ is \textbf{guided} by the pairs $(R_i, r_i)$.
 For a pointed space $X$ with a base point  $b$ the sequence   $\mathcal{F}$  is said to \textbf{respect
   the base point} if $\B(b, r_i)$ is contained in an element of ${\mathcal{F}}_{ij}$ for every  $j$.

\end{definition}
{\bf Remark.} For the results presented in this section and the following one it is
sufficient to assume $r_{i+1} > 2R_i$. The inequality $r_{i+1} > (100i+1)R_i$ will
be used in the last section for proving the factorization theorem.

\begin{definition}
Let $X$ be  a set. A sequence  $\mathcal{F}$ of  covers ${\mathcal{F}}_i, i \in \N,$ of    $X$ is said to \textbf{define
 $\asdim \leq n$} if   $\st {\mathcal{F}}_i$  refines ${\mathcal{F}}_{i+1}$,   ${\mathcal{F}}_{i}$ splits into the union of 
$n+1$ families ${\mathcal{F}}_{ij}, 1\leq j \leq n+1,$ of disjoint sets,  the sets of  ${\mathcal{F}}_{i+1, j}$
 are separated by ${\mathcal{F}}_i$ and for every $x,y \in X$ there is $i$ such that $y \in \st (x, {\mathcal{F}}_i )$.
 We say that $\mathcal{F}$ \textbf{separates the points} of a subset $Z$ of $X$ if ${\mathcal{F}}_1$ separates the points
 of $Z$.

\end{definition}
\begin{proposition}
\label{witnessing}
${}$

(i)
A metric space  of $\asdim \leq n$  admits a sequence of covers  witnessing $\asdim \leq n$.
A pointed metric space of $\asdim \leq n$ admits a sequence of covers witnessing $\asdim \leq n$
and respecting the base point.

(ii)
A sequence of covers witnessing $\asdim \leq n$ is also a sequence defining $\asdim \leq n$.
\end{proposition}
{\bf Proof.} 

(i) 
  follows from Proposition \ref{prop-for-witnessing}.
 
 (ii)  follows from the definitions of families witnessing $\asdim\leq n$ and defining $\asdim \leq n$.
 $\black$
\begin{proposition}
\label{metric}
Let ${\mathcal{F}}$ be  a sequence of covers ${\mathcal{F}}_i$ of a set $X$ with the splittings 
${\mathcal{F}}_{ij}$ such that $\mathcal{F}$  defines   $\asdim \leq n$. 
For $x, y \in X$ set $d_{\mathcal{F}}(x,y) $ to be the  maximal  $i$ such that ${\mathcal{F}}_i$ separates $x$ and $y$
if such $i$ exists and $d_{\mathcal{F}}(x,y)=0$ otherwise. Then for a subset   $Z\subset X$ separated by
${\mathcal{F}}$ we have:
${}$

(i) $d_{\mathcal{F}}$ is a metric on $Z$;

(ii)
${\mathcal{F}}_{i+1}$ restricted to $Z$ is $i$-bounded and  has a Lebesgue number  $i-1$, 
${\mathcal{F}}_{i+1, j}$ restricted to $Z$ is $(i-1)$-disjoint and,
as a result,  $\asdim Z \leq n$ (everything here with respect to  $d_{\mathcal{F}}$);

(iii)  if 
$\mathcal{F}$ is a sequence witnessing $\asdim \leq  n$ for a metric space $(X,d)$ then
$d$ and $d_{\mathcal{F}}$ are coarsely equivalent on $Z$

\end{proposition}
{\bf Proof.} 
${}$

(i) The only thing needed to be checked is the triangle inequality. Take $x,y, z \in Z$.
Clearly we may  assume that $i=d_{\mathcal{F}}(x,z)\geq d_{\mathcal{F}}(z, y) \geq 1$. Then
$x, y \in \st(z, {\mathcal{F}}_{i+1})$ and hence $x$ and $y$ are contained in an element of
${\mathcal{F}}_{i+2}$ and therefore $d_{\mathcal{F}}(x,y)\leq i+1$. Thus 
$d_{\mathcal{F}}(x,y)\leq  d_{\mathcal{F}}(x,z)+d_{\mathcal{F}}(z,y)$.

(ii) follows from the definition of $d_{\mathcal{F}}$ and  the definition of  a family defining $\asdim \leq n$.

(iii) Take a uniformly bounded cover $\mathcal{B}$  of $Z$ with respect to $d_{\mathcal{F}}$. By (ii), 
$\mathcal{B}$ refines ${\mathcal{F}}_i$ for some $i$ and therefore $\mathcal{B}$ is uniformly bounded with respect to $d$.
A similar argument also  works in the other direction.
 $\black$

\end{section}
\begin{section}{Universal space}
Let $(X^{\alpha}, d^\alpha)$ be a collection of  separable metric   spaces  with $\asdim \leq n$
representing, up to coarse equivalence, all the separable metric spaces with $\asdim \leq n$.  
 In each $X^\alpha$ we fix a base point
  $b^{\alpha}$ and,  by  (i) of Proposition  \ref{witnessing}, take  a sequence 
  ${\mathcal{F}}^\alpha$ of covers ${\mathcal{F}}_i^\alpha $ 
of $X^\alpha$  with the splittings ${\mathcal{F}}^\alpha_{ij}$ that
witnesses $\asdim \leq n$, respects   the base point $b^\alpha$
and
guided by the pairs $(R_i^\alpha, r_i^\alpha)$. Replacing  $X^\alpha$ by a coarsely equivalent
subset we may assume that ${\mathcal{F}}^\alpha$ separates the points of $X^\alpha$.

Denote by
 $X=\vee X^{\alpha}$  the wedge sum of the spaces $X^\alpha$ with the base point $b \in X$ obtained by identifying 
the base points of all $X^\alpha$.  We consider each $X^\alpha$ as a subset  of $X$
and denote by ${\mathcal{F}}$  the sequence of covers  $\mathcal{F}_i$ of  $X$  with the splittings ${\mathcal{F}}_{ij}$
defined as follows: ${\mathcal{F}}_{ij}$ is the union of ${\mathcal{F}}^\alpha_{ij}$ for all $\alpha$ 
where the sets containing the base point $b$
being replaced by their union.

\begin{proposition}
\label{wedge}
In the above setting the following holds:
${}$

(i) $\mathcal{F}$   is a sequence of covers  of $X$ defining $\asdim \leq  n$ and separating the points of $X$,
and hence, by (i)  and (ii)  of Proposition \ref{metric}, we have  $\asdim (X, d_{\mathcal{F}})\leq n$;

(ii)  each $(X^{\alpha}, d^{\alpha})$ is coarsely  embedded into $(X, d_{\mathcal{F}})$;

(iii) if a function $ f : X \lo Y$ to a metric space $Y$ isometrically embeds each  $X^\alpha$ into $Y$ with respect to
$d_{\mathcal{F}}$ then $f$ is coarsely continuous.

\end{proposition}
{\bf Proof}${}$

(i) Since ${\mathcal{F}}^\alpha$ separates the points of $X^\alpha$ we get that ${ \mathcal{F}}^\alpha_1$
consists of singletons.  Then ${\mathcal{F}}_1$ consists of singletons as well and therefore
${\mathcal{F}}$ separates the points of $X$.

 Let $x^\alpha \in X^\alpha$ and $x^\beta \in X^\beta$, and
let $i$ be such that $x^\alpha \in \B(b^\alpha, r^\beta_i)$ and $x^\beta \in \B(b^\beta, r_i^\beta)$.
Then $x^\alpha$ and $x^\beta$ are contained in an element of ${\mathcal{F}}_i$ because
${\mathcal{F}}^\alpha$ and ${\mathcal{F}}^\beta$ respect the base points.

Thus to show
 that ${\mathcal{F}}$ defines  $\asdim \leq n$ we only need to show that
 ${\mathcal{F}}_i$ refines ${\mathcal{F}}_{i+1}$ and separates ${\mathcal{F}}_{i+1, j}$.
Take  a point $x \in X$ and consider the following  cases.

 Case 1:  $b \notin \st(x, {\mathcal{F}}_i)$.  Then for $X^\alpha$ such that  $x \in X^\alpha$
we have  that the sets of ${\mathcal{F}}_i$  containing $x$ are exactly  the sets of ${\mathcal{F}}^\alpha_i$ 
containing $x$. By (ii) of Proposition \ref{witnessing},
  ${\mathcal{F}}^\alpha_i$ is also defining $\asdim \leq n$ and hence  $\st(x, {\mathcal{F}}^\alpha_i)$ is
  contained in an element  of ${\mathcal{F}}^\alpha_{i+1}$ and  no element of ${\mathcal{F}}^\alpha_i$ containing $x$ meets
  disjoint elements of ${\mathcal{F}}^\alpha_{i+1, j}$ for every $j$. Since
   ${\mathcal{F}}_i$ restricted to $X^\alpha$ coincides with ${\mathcal{F}}^\alpha_i$ we get  that
   $\st(x, {\mathcal{F}}_i)$ is
  contained in an element  of ${\mathcal{F}}_{i+1}$ and  no element of ${\mathcal{F}}_i$ containing $x$ meets
  disjoint elements of ${\mathcal{F}}_{i+1, j}$ for every $j$.

  Case 2:  $b \in \st(x, {\mathcal{F}}_i)$. Recall that ${\mathcal{F}}^\alpha$  witnesses $\asdim \leq n$.
  Then, since $r^\alpha_{i+1} > 2R^{\alpha}_i$, we have that  $\st(x, {\mathcal{F}}_i )$ is contained in the union
  of the balls $\B(b^\alpha, r^\alpha_{i+1})$ for all $\alpha$ and this union  in its turn is contained in an element of
  ${\mathcal{F}}_{i+1, j}$ for every $j$  because each ${\mathcal{F}}^\alpha$ respects the base point.
  Thus we get  that
   $\st(x, {\mathcal{F}}_i)$ is
  contained in an element  of ${\mathcal{F}}_{i+1}$ and  no element of ${\mathcal{F}}_i$ containing $x$ meets
  disjoint elements of ${\mathcal{F}}_{i+1, j}$ for every $j$.
\\

(ii) follows from (iii) of Proposition \ref{metric} and the fact that ${\mathcal{F}}_i$ restricted to $X^\alpha$
coincides with ${\mathcal{F}}^\alpha_i$.
\\

(iii) Consider   $F \in {\mathcal{F}}_i$. By (ii) of  Proposition \ref{metric}, ${\mathcal{F}}_i$ is $(i-1)$-bounded
with respect to $d_{\mathcal{F}}$.
  If $F$  does not contain  $b$ then $F$ is contained in some $X^\alpha$ and therefore
$\diam f(F)=\diam F \leq i-1$. If $F$ does contain $b$ then $F$ is  contained in
the union of the  elements 
of $ {\mathcal{F}}^\alpha_i$ containing $b$  for all $\alpha$ and therefore $\diam f(F) \leq 2(i-1)$.
   Thus $f({\mathcal{F}}_i)$ is uniformly   bounded. Moreover,  by (ii) of Proposition \ref{metric},
any uniformly bounded cover of $(X, d_{\mathcal{F}})$ refines ${\mathcal{F}}_i$ for some $i$ and, hence,
 $f$ is coarsely continuous.
 $\black$
 \\
 \\
 {\bf Proof of Theorem \ref{universal}.} Consider  the space $X$
 from Proposition \ref{wedge} and  assume that $X$ is equipped with the metric $d_{\mathcal{F}}$.
  Let $\U$ be  the Urysohn space.
 Recall that each  separable metric space isometrically embeds into $\U$
 and $\U$  is homogeneous  by isometries.  Then  one can define isometric embeddings $f^\alpha: X^\alpha \lo \U$
 (with respect to $d_{\mathcal{F}}$ restricted to $X^\alpha$)
 sending all the base points $b^\alpha$ to the same point in $\U$ and this way to define 
 a function $f : X \lo \U$
 that isometrically embeds each $X^\alpha \subset X$ into $\U$ (with respect to $d_{\mathcal{F}}$).
 By (iii) of Proposition \ref{wedge}, $f$ is coarsely continuous. Apply Theorem \ref{factorization}
 to factorize $f$ through coarsely continuous maps $g : X \lo Z$ and
 $h: Z \lo \U$ with $Z$ being separable metric with $\asdim  Z \leq \asdim X$. Recall that,
by (i) of Proposition \ref{wedge}, $\asdim X \leq n$ and hence $\asdim Z \leq n$.

Since $f$ coarsely (even isometrically)  embeds $X^\alpha$ into $\U$ we get that 
$g$ coarsely embeds $X^\alpha$ into $Z$. Indeed, if ${\mathcal{B}}$ is a uniformly bounded
cover of $Z$ then $h(\mathcal{B}) $ is uniformly bounded in $\U$ and, hence, the cover
$g^{-1}({\mathcal{B}})$ is a uniformly bounded on $X^\alpha$ since $g^{-1}({\mathcal{B}})$
 and $f^{-1}(h(\mathcal{B}))$ coincide on  $X^\alpha$   and 
$f^{-1}(h({\mathcal{B}}))$ is uniformly  bounded on $X^\alpha$ with respect to $d_{\mathcal{F}}$.
Finally note that, by (ii) of Proposition \ref{wedge}, the metrics $d^\alpha$ and
$d_{\mathcal{F}}$ are coarsely equivalent on  $X^\alpha$ and this shows that
$Z$ is a universal space for separable metric spaces of $\asdim \leq n$. 
$\black$ 
\end{section}

\begin{section}{Factorization Theorem}
Let us first make the following observation:
\begin{proposition}
\label{close}
Let $f : X \lo Y$, $g: X \lo Z$ and $h: Z \lo Y$ be coarsely continuous maps of metric spaces
such that  $f$ and $h \circ g$ are coarsely close  and $\w Y$ is infinite.  Then there is  a  metric  space
$Z' $ and coarsely continuous maps $g' : X \lo Z'$ and $h' : Z' \lo Y$ such that
$\asdim Z' \leq \asdim Z$, $\w Z' \leq \max \{  \w   Z, \w Y \}$ and $f=h' \circ g'$.

\end{proposition}
{\bf Proof.}  Set $Z'= \{ (f(x), g(x)) : x \in X \} \subset Y \times Z$ and
consider 
$Z'$ with the  metric  inherited from $Y \times Z$ and defined as
the maximum of the   distances  between the coordinates.
  Define $g' : X \lo Z'$
by $g'(x)=(f(x), g(x)), x \in X$ and let $h' : Z' \lo Y$ and $\pi: Z' \lo Z$ be the projections.
Clearly $g', h'$ and $\pi$ are coarsely continuous,
$f=h' \circ g'$ and $\w Z' \leq \max \{  \w   Z, \w Y \}$. 

Let us show that
$\pi$ is a coarse embedding. Take a uniformly bounded  cover 
$\mathcal{B}$ of $Z$. 
Then, since $h({\mathcal{B}})$ is
uniformly bounded and  $f$ and $h\circ g$ are coarsely close, we get
that $f(g^{-1}({\mathcal{B}}))$ is uniformly bounded as well. Note that
$\pi^{-1}(B)\subset f(g^{-1}(B)) \times B $  for every $B \in \mathcal{B}$ and 
hence
$\pi^{-1}({\mathcal{B}})$  is uniformly bounded. This implies that $\pi$ is
a coarse embedding and hence $\asdim Z' \leq \asdim Z$. $\black$
\\
\\
{\bf Proof of Theorem \ref{factorization}.} The theorem is trivial if $\w Y$ is finite, so we may 
assume that $\w Y$ is infinite. Let $\asdim X \leq n$. Fix  a base point $b$ in $X$. 
Take a sequence  ${\mathcal{F}}^X$ of  covers ${\mathcal{F}}^X_i$ of $X$ 
with splittings ${\mathcal{F}}^X_{ij}$ guided 
by  the pairs $(R_i, r_i)$ and witnessing $\asdim  \leq n$. 
For each $i$ take a uniformly bounded cover ${\mathcal{F}}^Y_i$ 
of $Y$  that
is refined by $f({\mathcal{F}}^X_i)$ with the cardinality of ${\mathcal{F}}^Y_i$ bounded by $\w Y$
  and  define ${\mathcal{F}}^0_{ij}$ as a  collection of disjoint subsets of $X$ such that each element 
  of ${\mathcal{F}}^0_{ij}$ is  a  union of elements
  of ${\mathcal{F}}^X_{ij}$ contained in an element of  $f^{-1}({\mathcal{F}}_i^Y)$,  each set of ${\mathcal{F}}^X_{ij}$ 
  is contained in some element of ${\mathcal{F}}^0_{ij}$ and  no element of ${\mathcal{F}}^X_{ij}$ is contained 
  in different elements of ${\mathcal{F}}^0_{ij}$.  
  Note that $f({\mathcal{F}}_{ij}^0)$ refines ${\mathcal{F}}^Y_i$.
  
  We will construct by induction collections ${\mathcal{F}}^p_{ij}, p \in \N,$ of subsets of $X$ satisfying the following conditions
  (the families ${\mathcal{F}}^p$ and ${\mathcal{F}}^p_i$  below are determined in the standard way by ${\mathcal{F}}^p_{ij}$, namely
  ${\mathcal{F}}^p_i$ is the union of the families
  ${\mathcal{F}}^p_{ij}$ for $1\leq j \leq n+1$ and 
  ${\mathcal{F}}^p$ is the sequence of the covers ${\mathcal{F}}^p_i$ of $X$):
  \\

  ($\dagger 1$) the cardinality of ${\mathcal{F}}^p_{ij}$ is bounded by  $\w Y$;
  
  ($\dagger 2$) the elements of ${\mathcal{F}}^p_{ij}$ are disjoint, each element of ${\mathcal{F}}^p_{ij}$ is a  union of 
  elements of ${\mathcal{F}}^X_{ij}$ and  each element of ${\mathcal{F}}^X_{ij}$ is
  contained in some element  of ${\mathcal{F}}^p_{ij}$;

  ($\dagger 3$) ${\mathcal{F}}^{p+1}_{ij}$ and ${\mathcal{F}}^p_{ij}$ restricted to $\B(b, ir_{i+p})$ coincide and
  ${\mathcal{F}}^{p+1}_{ij}$ refines ${\mathcal{F}}^p_{ij}$;

  ($\dagger 4$) $\st{\mathcal{F}}^{p+1}_{i}$ refines ${\mathcal{F}}^p_{i+1}$;

 ($\dagger 5$) the elements of ${\mathcal{F}}^p_{i+1,j}$ are separated  by ${\mathcal{F}}^{p+1}_{i}$.
  \\\\
  Clearly the relevant properties hold for ${\mathcal{F}}^0$. Assume that the construction is completed for 
  ${\mathcal{F}}^p_{i}$ with $p+i=m$ and proceed to $m+1$ in the following order
  ${\mathcal{F}}^0_{m+1}, {\mathcal{F}}^1_m, {\mathcal{F}}^2_{m-1}, \dots , {\mathcal{F}}^m_1$.
  Recall that ${\mathcal{F}}^0_{m+1}$ is already defined and assume that
  ${\mathcal{F}}^0_{m+1}, \dots,{\mathcal{F}}_{i+1}^{m-i}$ are already constructed.
  We will construct   ${\mathcal{F}}_{i}^{m-i+1}$ as follows.
  Denote  $p=m-i$.
  \\
  
  Let $1\leq j, t \leq n+1$ and
  let a set $A$ be the union of  some elements  of ${\mathcal{F}}^X_{ij}$.
  By splitting $A$ by   ${\mathcal{F}}^{p  }_{i+1,  t }$ we mean replacing $A$ by 
  the  family of disjoint subsets of $X$  which is the union of the following collections:
  \\
  
  collection 1: 
   the union of the elements of ${\mathcal{F}}^X_{ij}$
  contained in $A$   and intersecting $\B(b, ir_m)$  (collection 1 consists  of  only one set);
  
  collection 2:
  for each element  $F$ of ${\mathcal{F}}^{p}_{i+1, t}$
   take the union 
   of the elements of ${\mathcal{F}}^X_{ij}$
  which are  contained in both $A$ and  $F$ and  which were not used  in constructing collection 1
  (by ($\dagger 1$) the cardinality of collection 2 is bounded by $\w Y$);

  collection 3:
  for each element  $F$ of ${\mathcal{F}}^{p}_{i+1, t}$
   take the union 
   of the elements of ${\mathcal{F}}^X_{ij}$   which are
  contained  in  $A$ and  $(r_{i+1}/10)$-close  to $F$  and which  were not used  in constructing collections 1  and 2
   (by ($\dagger 1$) the cardinality of collection 3 is bounded by $\w Y$, and the elements of collection $3$ are disjoint because ($\dagger 2$) implies that ${\mathcal{F}}^{p}_{i+1,t}$ is $r_{i+1}$-disjoint);
  
  collection 4: the union of the elements of ${\mathcal{F}}^X_{ij}$
  contained in $A$  and that were not used in constructing collections 1, 2 and 3 (collection 4 consists  of  only one set).
  \\
  
  Now split the elements of ${\mathcal{F}}^p_{ij}$ by ${\mathcal{F}}^{p  }_{i+1,  1}$, the resulting family split by ${\mathcal{F}}^{p  }_{i+1,  2}$
  and proceed by induction  to the last splitting by ${\mathcal{F}}^{p  }_{i+1,  n+1}$  to finally get the family which we denote
  by ${\mathcal{F}}^{p+1}_{ij}$. The construction is completed. 
  
  \begin{claim}
  \label{conditions}
 All the  $\dagger$-conditions hold.
  \end{claim}
  
  We will prove Claim \ref{conditions} later. Let us show now how to derive from
  Claim \ref{conditions}  the proof of the factorization theorem.
   Condition ($\dagger 3$) implies that there is a unique collection ${\mathcal{F}}_{ij}$ of subsets
   of $X$ such that ${\mathcal{F}}_{ij}$ coincides with ${\mathcal{F}}^p_{ij}$ on the ball $\B(b, ir_{i+p})$ for
   every $p$.
   Condition ($\dagger 2$) implies that ${\mathcal{F}}_{ij}$ is a family of disjoint sets such that
   ${\mathcal{F}}_i={\mathcal{F}}_{i, 1} \cup \dots \cup {\mathcal{F}}_{i, n+1}$ covers $X$.
   Let $\mathcal{F}$ be the sequence of the covers  ${\mathcal{F}}_i$ with the splittings ${\mathcal{F}}_{ij}$.

 \begin{claim}
 \label{family}
 ${}$

 (i)  
${\mathcal{F}}^X_i$ refines ${\mathcal{F}}_i$
  and ${\mathcal{F}}_i$ refines ${\mathcal{F}}^p_i$  for every $p$;

  (ii)
 ${\mathcal{F}}$  defines  $\asdim \leq n$.
 
  \end{claim}
   {\bf Proof:} 
   
   (i) follows from ($\dagger 2$) and ($\dagger 3$).
   
   (ii)  Since ${\mathcal{F}}^X$ is a sequence witnessing  $\asdim \leq n$, any pair of points
   is contained in an element of ${\mathcal{F}}^X_i$ for some $i$. Then, by ($\dagger 1$) 
  this pair is  contained in an element  of   ${\mathcal{F}}^p_i$  for every $p$ and  hence,  by ($\dagger 3$),
  in an element of ${\mathcal{F}}_i$. Thus in order to show that ${\mathcal{F}}$ defines $\asdim \leq n$ we only
  need to show that
    $\st {\mathcal{F}}_i$ refines ${\mathcal{F}}_{i+1}$ and ${\mathcal{F}}_i$ separates ${\mathcal{F}}_{i+1,j}$.

   Take a point $x \in X$.
   By ($\dagger 3$)  and ($\dagger 4$), for every $p$  there is $F^p \in {\mathcal{F}}_{i+1}$ such that
   $\st(x, {\mathcal{F}}_i) \cap \B(b, ir_{i+p}) \subset  F^p \cap \B(b, ir_{i+p})$. Note that 
   $x$ belongs to at most $n+1$ elements of ${\mathcal{F}}_{i+1}$ and hence there is $F$ in 
   ${\mathcal{F}}_{i+1}$ such that $F=F^p$ for infinitely many $p$ and then
   $\st(x, {\mathcal{F}}_i)  \subset F$. Thus $\st {\mathcal{F}}_i$ refines ${\mathcal{F}}_{i+1}$.

   The property that ${\mathcal{F}}_i$ separates ${\mathcal{F}}_{i+1,j}$ follows from ($\dagger 3$)
   and  ($\dagger 5$).  
  $\black$
  \\
  
  Let $Z$ be a maximal  subset of $X$ separated  by ${\mathcal{F}}_1$. Define
   a function $g: X \lo Z$ by sending $x \in X$ to a point $z \in Z$  such that $x \in \st(z, {\mathcal{F}}_1)$.
   Consider $Z$ with the metric $d_{\mathcal{F}}$  determined by ${\mathcal{F}}$ as described 
  in  Proposition \ref{metric}.

\begin{claim} 
\label{last} We have that $\w Z \leq \w Y$,
$\asdim Z \leq n$, the functions $g: X \lo Z$ and   $h=f|Z  : Z \lo Y$ are coarsely continuous
(everything here with respect to $d_{\mathcal{F}}$),  and  $f$  and    $h\circ g$ are coarsely close.

\end{claim}
  {\bf Proof.}  Since ${\mathcal{F}}_1$ separates  the points of $Z$,  we get that ${\mathcal{F}}_1$ restricted to $Z$ consists of
  singletons and therefore the cardinality of $Z$ is bounded by the cardinality of ${\mathcal{F}}_1$.
  Then, since ${\mathcal{F}}_1$ restricted  to $\B(b, r_{1+p})$  coincides 
  with ${\mathcal{F}}^p_1$ for every $p$, we get by  (${\dagger 1}$)  that the cardinality of ${\mathcal{F}}_1$ is bounded by $\w Y$ and hence
   $\w Z \leq \w Y$.
  The property  $\asdim \leq n$ follows from (ii) of Proposition \ref{metric}. 
  
  Let us show that $g$ is coarsely continuous.
  Take $F\in {\mathcal{F}}^X_i$ and $x_1, x_2 \in F$ and let $z_1=g(x_1)$ and $z_2=g(x_2)$.
  Then $x_1 \in \st(z_1, {\mathcal{F}}_1)$ and $x_2 \in \st(z_2, {\mathcal{F}}_1)$. By (ii) of Claim \ref{family},
 $ x_1 \in \st(z_1, {\mathcal{F}}_i)$ and $x_2 \in \st(z_2, {\mathcal{F}}_i)$ and  then
 $F \subset  \st(z_1, {\mathcal{F}}_{i+1})$  and  $F \subset  \st(z_2, {\mathcal{F}}_{i+1})$. Therefore,  by (ii) of Proposition \ref{metric},
 $z_1$ and $z_2$ are $(i+1)$-close with respect to $d_{\mathcal{F}}$. Thus,    $g({\mathcal{F}}^X_i)$ is  uniformly bounded
 with respect to $d_{\mathcal{F}}$ and hence  $g$ is coarsely continuous.
 
    Let us show that $h$ is coarsely continuous. Recall that $f({\mathcal{F}}^0_i)$ refines ${\mathcal{F}}^Y_i$. Then,
    by (i) of Claim \ref{family}, ${\mathcal{F}}_i$ refines ${\mathcal{F}}^0_i$ and we get that  $f({\mathcal{F}}_i)$ is uniformly bounded
    and hence, again by (ii) of Proposition \ref{metric}, $h$ is coarsely continuous.
    
    Let us finally show that $f$ and $h \circ g$ are coarsely close. Take a point $x \in X$ and let $z=g(x)$.
    Since $x \in \st (z,  {\mathcal{F}}_1)$, ${\mathcal{F}}_1$ refines ${\mathcal{F}}^0_1$ and $f({\mathcal{F}}^0_1)$
    refines ${\mathcal{F}}^Y_1$ we get that $f(x)$ and $h(g(x))=f(z)$ are contained in an element of ${\mathcal{F}}^Y_1$
    and hence $f$ and $h\circ g$ are coarsely close. $\black$
 \\
 \\
  Thus Theorem \ref{factorization} follows from  Claim  \ref {last} and Proposition \ref{close}. The only thing left is:
  \\
  \\
  {\bf Proof  of Claim \ref{conditions}}.
   Conditions ($\dagger 1$), ($\dagger 2$) and ($\dagger 3$) obviously follow from  the construction. So the only conditions
  we need to verify are ($\dagger 4$) and ($\dagger 5$) for ${\mathcal{F}}^{p+1}_i$ assuming that  all the $\dagger$-conditions hold
  for all the covers constructed before ${\mathcal{F}}^{p+1}_i$. Recall that $m=p+i$.
  Note that, by ($\dagger 2$), 
  \\
  
  (*) $r_{i+1}$ is a Lebesgue number of ${\mathcal{F}}^p_{i+1}$ and

 (**)  ${\mathcal{F}}^p_{i+1,t}$ is $r_{i+1}$-disjoint  for every $t$.
 \\
Also recall that

 (***) ${\mathcal{F}}^X_i$ is
  $R_i$-bounded with $ (100i+1)R_i< r_{i+1}$ (in particular $R_i < r_{1+1}/10$).
  \\
  
  Fix a point $x \in X$. We will say
  that $x$ satisfies ($\dagger 4$) if $\st( x, {\mathcal{F}}^{p+1}_i)$ is contained  in an element
  of ${\mathcal{F}}^p_{i+1}$ and we will say that $x$ satisfies ($\dagger 5$) if for every $j$ no
  element of ${\mathcal{F}}^{p+1}_i$ containing $x$  meets disjoint elements of ${\mathcal{F}}^p_{i+1,j}$
  one of which contains $x$. Note that if every point of $X$ satisfies ($\dagger 4$) and ($\dagger 5$)
  then ${\mathcal{F}}^{p+1}_i$ satisfies ($\dagger 4$) and ($\dagger 5$) as well.
  
   Let $E$ be an element of ${\mathcal{F}}^{p+1}_i$.
   Recall that $E$ is constructed from an element of ${\mathcal{F}}^p_i$  by  a  sequence  of splittings
    by the families ${\mathcal{F}}^p_{i, 1}, \dots {\mathcal{F}}^p_{i, n+1}$.
   On each  step of this construction we create  4 collections (collections 1, 2, 3 and 4). Only one
   of these collections has an element containing $E$ and  we will refer  to    this collection as the collection
   refined by $E$. We will refer to the step of splitting by ${\mathcal{F}}^p_{i, t}$ in the construction of $E$ as 
   step $t$. 
   For  given $x$ and $t$ we assume throughout  the proof  that 
\\

 ($\diamond$) \ \  $E$ is an element of $\mathcal{F}^{p+1}_{i}$ that contains $x$, $\mathcal{E}$ is the collection refined by $E$ created on step $t$
 of the construction of $E$, and 
 $E_{\mathcal{E}}$ is the unique element of $\mathcal{E}$ that contains $E$.
  \\\\

   Consider the following cases.
  \\
  
  {\bf Case 1:} $\st( x, {\mathcal{F}}^X_i)$  does not meet $\B(b, ir_m)$.
  \\
  
  We will  show that $x$ satisfies ($\dagger 4$).
  By (*) and (***),  take an element $F \in {\mathcal{F}}^p_{i+1, t}$ containing $\st( x, {\mathcal{F}}^X_i)$. 
  Having $x$ and $t$ chosen we let $E$, $\mathcal{E}$, and $E_{\mathcal{E}}$ be as in $(\diamond)$.
   By the assumption of  Case 1, $\mathcal{E}$ cannot be collection 1.
  Then ${\mathcal{E}}$ must be collection 2  with $E_{\mathcal{E}} \subset F$.
   Thus $E \subset F$ and hence 
 $x$ satisfies ($\dagger 4$).
  \\
  
  Now we will show that $x$ satisfies ($\dagger 5$).
   Fix any $t$  and  suppose $F$   is an element of ${\mathcal{F}}^p_{i+1, t}$ containing
  $x$.
   Having $x$ and $t$ chosen we let $E$, $\mathcal{E}$, and $E_{\mathcal{E}}$ be as in $(\diamond)$.
   By the assumption of Case 1, ${\mathcal{E}}$ is not collection 1.
   Then, by (**) and (***), ${\mathcal{E}}$ must be either collection 2 or collection 3 and
   then $E_{\mathcal{E}}$ is either contained in $F$ or $(r_{i+1}/10)$-close to $F$. Hence, by (**),
   $E_{\mathcal{E}}$ meets no element of ${\mathcal{F}}^p_{i+1, t}$ different from $F$. 
 Thus $x$ satisfies ($\dagger 5$).
  \\
  \\

 {\bf  Case 2:} $\st( x, {\mathcal{F}}^X_i)$ does   meet $\B(b, ir_m)$ and $p=0$. Note that in this case $m=i$ and
  hence $2(ir_m + 2R_i )< r_{i+1}/10$.
  \\
  
  We will show that $x$ satisfies ($\dagger 4$).  
   By (*) take an element $F \in {\mathcal{F}}^p_{i+1,t}$ containing  $ \B(b, ir_m +2R_i)$.
  Then, by (***), $\st(x, {\mathcal{F}}^X_i)\subset F$.
     Having $x$ and $t$ chosen we let $E$, $\mathcal{E}$, and $E_{\mathcal{E}}$ be as in $(\diamond)$.
   If ${\mathcal{E}}$ is collection 1 then, again by (***),  $ E_{\mathcal{E}} \subset \B(b, ir_m +2R_i)$ and 
   hence  $ E_{\mathcal{E}} \subset F$.  If $\mathcal{E}$ is not collection 1
  then, 
   since $\st(x, {\mathcal{F}}^X_i)\subset F$, the only option left is that $\mathcal{E}$ is collection 2 and
  hence  $ E_{\mathcal{E}} \subset F$.
Thus $x$ satisfies ($\dagger 4$).
  \\
  
  Now let us show that $x$ satisfies ($\dagger 5$).
   Fix any $t$  and  suppose  $F$   is an element of ${\mathcal{F}}^p_{i+1, t}$ containing
   $x$.   
    Having $x$ and $t$ chosen we let $E$, $\mathcal{E}$, and $E_{\mathcal{E}}$ be as in $(\diamond)$.
  If $\mathcal{E}$ is collection 1 then 
   $E_{\mathcal{E}} $ is contained in $\B(b, ir_m +2R_i)$ and hence $E_{\mathcal{E}}$ is $(r_{i+1}/10)$-close to $F$.
   If  $\mathcal{E}$ is either collection 2 or collection 3 and then
   in both cases  $E_{\mathcal{E}}$ is $(r_{i+1}/10)$-close to $F$. By (**) and (***), $\mathcal{E}$ cannot be collection 4.
   Thus, by (**),  $E_{\mathcal{E}}$ cannot meet an element of ${\mathcal{F}}^p_{i+1, t}$ different from $F$
   and therefore $x$ satisfies ($\dagger 5$).
   \\
   
   {\bf  Case 3:} $\st( x, {\mathcal{F}}^X_i)$ does   meet $\B(b, ir_m)$ and $p>0$. 
  Note that  then, by (***),  $\st( x, {\mathcal{F}}^X_i) \subset \B(b, (i+1)r_m)$. 
  \\
  
  Let us show that $x$ satisfies ($\dagger 4$). By the inductive assumption for ($\dagger 4$),
  $\st {\mathcal{F}}^p_i$ refines ${\mathcal{F}}^{p-1}_{i+1}$. Then,
  by   ($\dagger 3$), there  is $F\in {\mathcal{F}}^{p}_{i+1}$
  such that 
  \\
  
 ($\bullet$) \ \  $\st(x, {\mathcal{F}}^p_i) \cap \B(b,( i+1)r_m ) \subset  F \cap \B(b,( i+1)r_m ) $.
  \\\\
  Note that then, by  (***) and ($\dagger 2$), we also have  $\st(x, {\mathcal{F}}^X_i) \subset F$.
  Suppose  $F \in {\mathcal{F}}^p_{i+1, t}$.
     Having $x$ and $t$ chosen we let $E$, $\mathcal{E}$, and $E_{\mathcal{E}}$ be as in $(\diamond)$. If ${\mathcal{E}}$ is collection 1  then
 $E_{\mathcal{E}} \subset \B(b, (i+1)r_m)$ and, then 
   by ($\dagger 3$),  $E$ is  contained in an element of ${\mathcal{F}}^p_i $, and
   hence, by ($\bullet$),  $ E \subset F$. If $\mathcal{E}$ is not collection 1 then, since $\st(x, {\mathcal{F}}^X_i) \subset F$,
    the only option left
   is that ${\mathcal{E}}$ is collection 2  and hence $E_{\mathcal{E}} \subset F$.
   Thus $x$ satisfies ($\dagger 4$).
   \\
   
   Now we will show that $x$  satisfies ($\dagger 5$). Fix any $t$ and suppose  $F$ is an element of  ${\mathcal{F}}^p_{i+1, t}$
   containing $x$. 
      Having $x$ and $t$ chosen we let $E$, $\mathcal{E}$, and $E_{\mathcal{E}}$ be as in $(\diamond)$.
   If $\mathcal{E}$ is collection 1 then
     $E_{\mathcal{E}} \subset \B(b, (i+1)r_m)$  and,  hence by ($\dagger 3$),  $E$ is contained in an element of ${\mathcal{F}}^p_i$.
     Thus, by the inductive assumption for  ($\dagger 5$),  $E$ 
     cannot meet disjoint
   elements of ${\mathcal{F}}^{p-1}_{i+1, t}$ and hence, by ($\dagger 3$),   $E$ cannot meet disjoint elements of
    ${\mathcal{F}}^{p}_{i+1, t}$  as well.
   If ${\mathcal{E}}$ is collection 2 or collection 3 
   then $E_{\mathcal{E}}$ is $(r_{i+1}/10)$-close to $F$ and, by (**), $E_{\mathcal{E}}$  cannot
   meet disjoint elements of  ${\mathcal{F}}^{p}_{i+1, t}$.
   By (**) and (***), $\mathcal{E}$ cannot be collection 4.
   Thus $F$ is the only element
   ${\mathcal{F}}^p_{i+1, t}$ that $E$ meets  and therefore
    $x$ satisfies ($\dagger 5$).\\
    
    Thus  every point of $X$ satisfies ($\dagger 4$) and ($\dagger 5$) and
   the claim is proved. $\black$

  \end{section}

 Jerzy Dydak\\
University of Tennessee, Knoxville, TN 37996, USA and\\
Xi'an Technological University, No.2 Xuefu zhong lu, Weiyang district, Xi'an, China 710021
\\
jdydak@utk.edu
\\
\\
Michael Levin\\
Ben Gurion University of the Negev, Beer Sheva, 8410501, Israel\\
mlevine@math.bgu.ac.il
\\
\\
Jeremy Siegert\\
Ben Gurion University of the Negev, Beer Sheva, 8410501, Israel\\
siegertj@post.bgu.ac.il

\end{document}